\documentclass[preprint,12pt]{elsarticle}
\usepackage{float,epsfig,bm,floatflt,here}
\usepackage{rotating}
\usepackage{amssymb}
\usepackage{amsmath}

\journal{Journal of Computational Physics}

\begin{document}

\begin{frontmatter}
\title{\textbf{Fourth order compact schemes for variable coefficient parabolic problems
with mixed derivatives}}

\author{Shuvam Sen}
\ead{shuvam@tezu.ernet.in}
\address{Department of Mathematical Sciences, Tezpur University,\\ Tezpur 784028, Assam, INDIA}

\begin{abstract}
In this article, we have developed a higher order compact numerical method for variable coefficient parabolic problems with mixed derivatives. The finite difference scheme, presented here for two-dimensional domains, is based on fourth order spatial discretization. The time discretization has been carried out using using second order Crank-Nicolson. The present scheme shows good dispersion relation preserving property and has been thoroughly investigated for stability. The discrete Fourier analysis shows that the method is unconditionally stable. The fact that the method has been particularly developed for parabolic equations with mixed derivatives makes it suitable for solving incompressible Navier-Stokes (N-S) equations in irregular domains. To verify the proposed method, several problems with exact and benchmark solutions has been investigated. The proposed compact discretization has been extended to tackle flows of varying complexities governed by the two-dimensional unsteady N-S equations in domain beyond rectangular. The results show good agreement for all the problems considered.
\end{abstract}


\end{frontmatter}

\section{Introduction}
\subsection{Problem formulation}
Let $\Omega$ be a rectangular domain in $\mathbb{R}^2$ and $(0, T]$ be the time interval with $T>0$. The paper describes a compact fourth order numerical method for solving parabolic partial differential equation (PDE) of the form
\begin{eqnarray}\label{1}
\begin{cases}
  \partial_t\phi(X,t)+A\phi(X,t)=s(X,t),\;\;\;\;\;\;\;\;\;\;(X,t)\in\Omega\times(0,T] \\
  \phi(X,0)=\phi_0(X),\;\;\;\;\;\;\;\;\;\;\;\;\;\;\;\;\;\;\;\;\;\;\;\;\;\;\;\;\;\;\;X\in\Omega  \\
  b_1(X,t)\phi+b_2(X,t)\partial_n\phi=g(X,t),\;\;\;\;\;\;X\in\partial\Omega,\;\;t\in(0,T]
\end{cases}
\end{eqnarray}
for the unknown transport variable $\phi(X,t)$ defined over $\Omega\times(0,T]\subset\mathbb{R}^2\times\mathbb{R}$. Here $A$ is a variable coefficient partial differential operator defined as
\begin{eqnarray}\label{2}
[A\phi](X,t)&=&-\alpha_1(X,t)[\partial_{xx}]\phi(X,t)-\beta(X,t)[\partial_{xy}]\phi(X,t)-\alpha_2(X,t)[\partial_{yy}]\phi(X,t)\nonumber\\
&&+c_1(X,t)[\partial_{x}]\phi(X,t)+c_2(X,t)[\partial_{y}]\phi(X,t)+d(X,t)\phi(X,t)
\end{eqnarray}
with $X=(x,y)$. Further the coefficients $\alpha_1(X,t)$, $\alpha_2(X,t)$, $\beta(X,t)$, $c_1(X,t)$, $c_2(X,t)$, $d(X,t)$ and forcing function $s(X,t)$ together with $\phi_0(X)$ and $g(X,t)$ are assumed to be sufficiently smooth. Here the only additional restriction we have on the Eq. (\ref{1}) is the positive definiteness of the diffusion matrix. This is equivalent to $\alpha_1(X,t)>0$, $\alpha_2(X,t)>0$ and $|\beta(X,t)|^2<4\alpha_1(X,t)\alpha_2(X,t)$ $\forall$ $(X,t)\in\Omega\times(0,T]$. $b_1$ and $b_2$ are arbitrary coefficients describing the boundary condition as a Dirichlet, Neumann, or Robin type in the boundary normal direction $n$.

The generalized convection-diffusion equations with mixed derivatives, given in Eq. (\ref{1}), arise in many applications. For example Heston equation which financial mathematician use for option pricing in stochastic volatility models \cite{hou_fou_10,dur_fou_12}. We also note the occurrence of such differential equations in mathematical biology \cite{hou_wel_07} and also in numerical mathematics when coordinate transformations are applied to convection-diffusion equations on non-rectangular domains \cite{mck_wal_wil_96,pan_kal_dal_07}. Such transformations allow us to work on simple rectangular domains or uniform grids although the original problem is considered on complex domains or with non-uniform grids which may be necessary to account for the stiff behaviour of solution in some part of the domain.

\subsection{Prior work}
In the mathematical literature there exist various different approaches to approximate solutions of PDE. The most popular one amongst these approaches, that has been historically used by the computational fluid dynamics (CFD) community, is the finite difference  (FD) method. A higher order compact (HOC) finite difference scheme is one which in two-dimensions employ a nine-point computational stencil using the eight directly adjacent nodes of the reference grid point and offer an accuracy of order four or higher. These schemes lead to a system of equations resulting in a coefficient matrix with smaller bandwidth as compared to non-compact schemes. Apart from solving convection-diffusion equation, different compact schemes have been used successfully to solve non-linear incompressible Navier-Stokes (N-S) equations in all three forms \emph{viz.} the stream function-vorticity \cite{pan_kal_dal_07,spo_car_95,li_tan_for_95,kal_dal_das_02,kal_das_dal_04,wan_zho_zha_06,ge_cao_11,sen_13}, the primitive variables \cite{kal_sen_07,tia_lia_yu_11} and the biharmonic \cite{ben_cro_fis_05,gup_kal_05,tia_yu_11,sen_kal_gup_13} formulations.

Here it is worthwhile to point out that although a plethora of compact schemes \cite{spo_car_95,kal_dal_das_02,kal_das_dal_04,wan_zho_zha_06,ge_cao_11,sen_13,gup_man_ste_84,spo_car_01,kar_zha_04,you_06,tia_ge_07,tia_11} have been developed for convection-diffusion equation but few can tackle the generalized one as specified in Eq. (\ref{1}). Of course a handful of contributions having linkage to the equation can be found in the works of Fourni\'{e} and Karaa \cite{fou_kar_06}, Pandit et al. \cite{pan_kal_dal_07}, Karaa \cite{kar_07} and D\"{u}ring and Fourni\'{e} \cite{dur_fou_12}. Fourni\'{e} and Karaa \cite{fou_kar_06} in 2006 derived a fourth order compact FD scheme for a two-dimensional (2D) elliptic PDE with mixed derivative by considering the PDE itself as an auxiliary relation. But in their work they restricted $\alpha_1=1=\alpha_2$ and also considered $\beta$ to be a constant with $\beta^2<4$. Application of this approach to more general problems such as Eq. (\ref{1}) is not straightforward and may not be possible. This is accentuated in the work of Pandit et al. \cite{pan_kal_dal_07} in 2007, where the authors tried to derive a compact approximation of parabolic equation. Here the authors were ultimately constrained to work with situations where the mixed derivative is absent. Karaa \cite{kar_07}, also in 2007, proposed a fourth order compact FD scheme for solving 2D elliptic and parabolic equations with mixed derivative having variable coefficient by using polynomial approximation, but was again limited by the choice of $\alpha_1=1=\alpha_2$. In 2012, D\"{u}ring and Fourni\'{e} \cite{dur_fou_12} used a compact scheme, having fourth order accuracy in space and second order accuracy in time, for option pricing using Heston model. In that manuscript the authors using a variable transformation arrived at a system with $\alpha_1=\alpha_2$. Note that such transformations are not always certain for a parabolic PDE with variable coefficients.

Although HOC approximation of Eq. (\ref{1}) is yet to be established, different splitting schemes and their alternating direction implicit (ADI) implementation can be found in the literature. Among them the works of in't Hout and Foulon \cite{hou_fou_10}, in't Hout and Welfert \cite{hou_wel_07}, Mckee et al. \cite{mck_wal_wil_96} and the references therein deserves special mention. Recently Martinsson \cite{mar_13} has designed a composite spectral collocation scheme for the equation with smooth solutions. These works clearly establishes importance of generalized parabolic equations. Here we will like to highlight that in many applications it is desirable to use higher order numerical methods to obtain accurate solution. Of late Sen \cite{sen_13} has developed a new family of implicit HOC schemes for unsteady convection-diffusion equation with variable convection coefficient. The schemes, where transport variable and its first derivatives are carried as the unknowns, combine virtues of compact discretization and Pad\'{e} approximation. In this manuscript we generalize this philosophy to propose a HOC formulation for variable coefficient parabolic problem with mixed derivative. We are able to obtain a new compact stable scheme with truncation error of order four in space and two in time. To the best of our knowledge compact schemes for generalized parabolic 2D convection-diffusion equations is not available in the literature and in this paper we intend to address the same. The schemes thus developed have been tested for their stability by using discrete Fourier analysis and shows better phase and amplitude error properties. We also plan to augment the scheme thus developed to solve incompressible N-S equations in irregular domains.

\subsection{Outline of this paper}
The paper is organized as follows. In the first part, we present the HOC difference scheme for Eq. (\ref{1}) and its extension to the 2D incompressible Navier-Stokes equations. In the second part, we apply this method to numerical computation to validate our approach. We consider some boundary layer and local singularity problems, as well as the classical impulsively started flow past circular problem as numerical examples to demonstrate the accuracy, effectiveness and efficiency of the present method. Finally, concluding remarks have been presented.

\section{Fourth order compact schemes for parabolic problem}
\subsection{Spatial compact discretization}
We begin by briefly discussing the development of HOC formulation for the steady state form of equation (\ref{1}), which is obtained when $\alpha_1$, $\beta$, $\alpha_2$, $c_1$, $c_2$, $d$, $s$ and $\phi$ are independent of $t$. Under these conditions, equation (\ref{1}) becomes
\begin{eqnarray}\label{3}
\begin{cases}
  A\phi(X)=s(X),\;\;\;\;\;\;\;\;\;\;\;\;\;\;\;\;\;\;\;\;\;\;\;\;\;X\in\Omega \\
  b_1(X)\phi+b_2(X)\partial_n\phi=g(X),\;\;\;\;\;\;X\in\partial\Omega
\end{cases}
\end{eqnarray}
For simplicity we assume $\Omega=[a_1,a_2]\times[a_3,a_4]$. In order to obtain a compact spatially fourth order accurate discretization we lay out a grid $a_1=x_0<x_1<...<x_M=a_2$, $a_3=y_0<y_1<...<y_N=a_4$ with $x_i=x_0+ih$ for $0\le i\le M$ and $y_j=y_0+jk$ for $0\le j\le N$. Consider the following approximations for second order space derivatives appearing in equation (\ref{3})
\begin{equation}\label{4}
\partial_{xx}\phi_{i,j}=2\delta^2_x\phi_{_{i,j}}-\delta_x\phi_{x_{i,j}}+O(h^4),
\end{equation}
\begin{equation}\label{5}
\partial_{yy}\phi_{i,j}=2\delta^2_y\phi_{_{i,j}}-\delta_y\phi_{y_{i,j}}+O(k^4),
\end{equation}
\begin{equation}\label{6}
\partial_{xy}\phi_{i,j}=\delta_{x}\phi_{y_{i,j}}+\delta_{y}\phi_{x_{i,j}}-\delta_{x}\delta_{y}\phi_{i,j}+O(h^2k^2).
\end{equation}
Here $\delta_x$, $\delta_y$, $\delta^2_{x}$ and $\delta^2_{y}$ are usual central difference operators and $\phi_{i,j}$ denote the approximate value of $\phi(X_{i,j})$ at a typical grid point $X_{i,j}=(x_i,y_j)$. Detailed derivation of the above discretizations can be found in \cite{sen_13}. We thus obtain an $O(h^4, k^4, h^2k^2)$ approximation for equation (\ref{3}) on a nine point stencil as
\begin{equation}\label{7}
A_{h,k}\phi_{i,j}=s_{i,j}
\end{equation}
where the discrete operator $A_{h,k}$ is defined as
\begin{eqnarray}\label{8}
A_{h,k}\phi_{i,j}&=&(-2\alpha_{1_{i,j}}\delta^2_x-2\alpha_{2_{i,j}}\delta^2_y+\beta_{_{i,j}}\delta_{x}\delta_{y}+d_{_{i,j}})\phi_{_{i,j}}\nonumber\\
&&+(\alpha_{1_{i,j}}\delta_x-\beta_{_{i,j}}\delta_y+c_{1_{i,j}})\phi_{x_{i,j}}+(\alpha_{2_{i,j}}\delta_y-\beta_{_{i,j}}\delta_x+c_{2_{i,j}})\phi_{y_{i,j}}.
\end{eqnarray}
The finite difference operator given above depends on the three grid functions $\phi$, $\phi_x$ and $\phi_y$. If we are interested in an approximation depending only on $\phi$, then we need compatible fourth order approximations for space derivatives $\phi_{x_{i,j}}$ and $\phi_{y_{i,j}}$ in terms of $\phi_{i,j}$. This is accomplished by using Pad\'{e} approximations
\begin{equation}\label{9}
\bigg(I+\frac{h^2}{6}\delta^2_{x}\bigg)\phi_{x_{i,j}}=\delta_{x}\phi_{i,j},
\end{equation}
and
\begin{equation}\label{10}
\bigg(I+\frac{k^2}{6}\delta^2_{y}\bigg)\phi_{y_{i,j}}=\delta_{y}\phi_{i,j}.
\end{equation}
Note that compared to standard HOC formulation \cite{kal_dal_das_02}, we are not required to approximate the derivatives of the convection coefficients $c_1$, $c_2$ and forcing function $s$. The equation (\ref{7}) can be viewed as a banded system with only nine non zero diagonals; of course drawback of requiring to approximate $\phi_{x_{i,j}}$ and $\phi_{y_{i,j}}$ separately using (\ref{9}) and (\ref{10}) respectively remain.

\subsection{Modified wave number analysis}
\begin{figure}[!p]
\begin{minipage}[b]{.6\linewidth}   \
\hfill
\centering\psfig{file=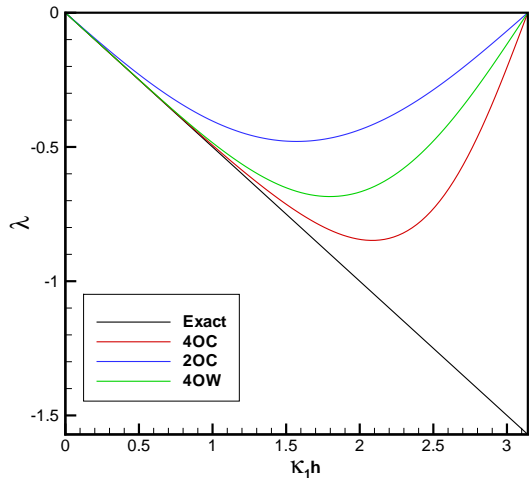,width=1.0\linewidth}
 (a)
\end{minipage}            \hspace{-2.5mm}
\begin{minipage}[b]{.6\linewidth}
\centering\psfig{file=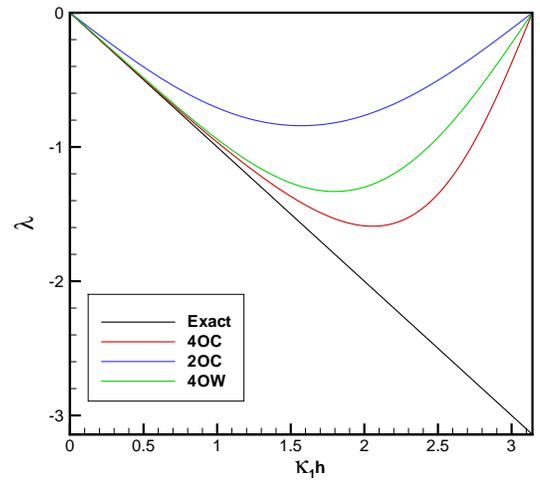,width=1.0\linewidth}
 (b)
\end{minipage}            \hspace{-2.5mm}
\vspace{.4cm}
 \begin{minipage}[b]{.6\linewidth}   \
\centering\psfig{file=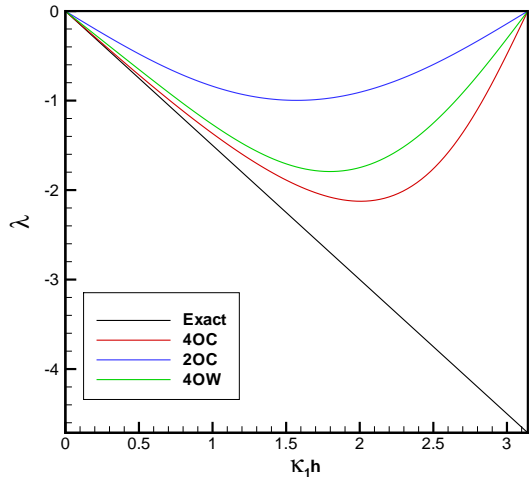,width=1.0\linewidth}
 (c)
\end{minipage}            \hspace{-2.5mm}
\begin{minipage}[b]{.6\linewidth}
\centering\psfig{file=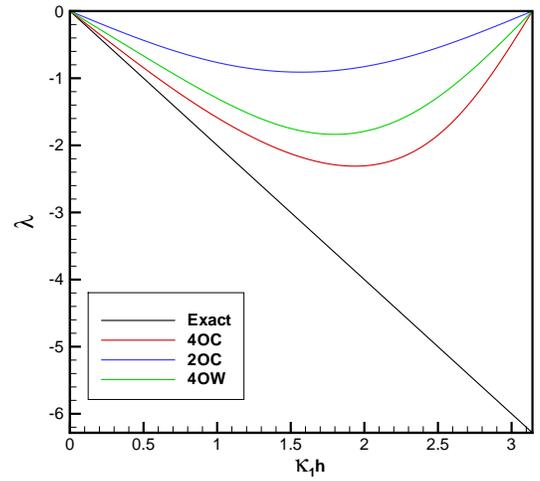,width=1.0\linewidth}
 (d)
\end{minipage}            \hspace{-2.5mm}
\begin{center}
\caption{{\sl Comparison of nondimensional characteristics for three different schemes with the exact one at (a) $\kappa_2k=0.5$, (b) $\kappa_2k=1.0$, (c) $\kappa_2k=1.5$, (d) $\kappa_2k=2.0$.} }
\label{fig:char}
\end{center}
\end{figure}

A detailed wave number analysis of the fourth order compact approximation for $\phi_{xx}$ was carried out by Sen \cite{sen_13}. Here we will like to examine the characteristic of the newly used fourth order compact approximation for the mixed derivative $\phi_{xy}$. We begin by considering the test equation
\begin{equation}\label{10.1}
\phi_{xy}=s(x,y).
\end{equation}
Consider the trial function $\tilde{\psi}=e^{I(\kappa_1 x+\kappa_2 y)}$ $(I=\sqrt{-1})$ where $\kappa_1$ and $\kappa_2$ are the wave numbers corresponding to $x$ and $y$ directions respectively. Then the exact characteristic of the equation (\ref{10.1}) is
\begin{equation}\label{10.2}
\lambda_{Exact}=-\kappa_1\kappa_2.
\end{equation}
Using the compact approximation given in Eq. (\ref{6}) in Eq. (\ref{10.1}) we get the corresponding discretized equation as
\begin{equation}\label{10.3}
\delta_{x}\phi_{y_{i,j}}+\delta_{y}\phi_{x_{i,j}}-\delta_{x}\delta_{y}\phi_{i,j}=s_{i,j}
\end{equation}
together with the Pad\'{e} approximation given in Eqs. (\ref{9}) and (\ref{10}). It is easy to see that the above fourth order compact discretization for mixed derivative has characteristic
\begin{equation}\label{10.4}
\lambda_{4OC-M}=-\frac{\sin (\kappa_1h)\sin (\kappa_2k)}{hk}\bigg[\frac{3}{2+\cos(\kappa_1h)}+\frac{3}{2+\cos(\kappa_2k)}-1\bigg].
\end{equation}
Note that the characteristics for the second order accurate central difference approximation and fourth order accurate wide stencil approximation for the mixed derivatives are
\begin{equation}\label{10.5}
\lambda_{2OC}=-\frac{\sin (\kappa_1h)\sin (\kappa_2k)}{hk}
\end{equation}
and
\begin{equation}\label{10.6}
\lambda_{4OW}=-\frac{\sin (\kappa_1h)\sin (\kappa_2k)}{9hk}(4-\cos(\kappa_1h))(4-\cos(\kappa_2k))
\end{equation}
respectively. To the best of our knowledge no HOC approximation for the mixed derivative is available in litarature. To get a clear idea of the dissipation errors associated with each of the above discretizations the non-dimensional characteristics as a function of $\kappa_1h$ corresponding to four different values of $\kappa_2k=0.5,1.0,1.5,2.0$ have been shown in figure \ref{fig:char}. Note that the expressions for characteristics being symmetrical with respect to $\kappa_1h$ and $\kappa_2k$ the above set of values give comprehensive idea of associated errors. The figure clearly indicates that the fourth order compact discretization discussed here has superior wave resolution property than the other discretization procedures that can be used for mixed derivatives.
\subsection{Implicit time discretization}
\begin{figure}[!h]
\begin{center}
\includegraphics[width=5in]{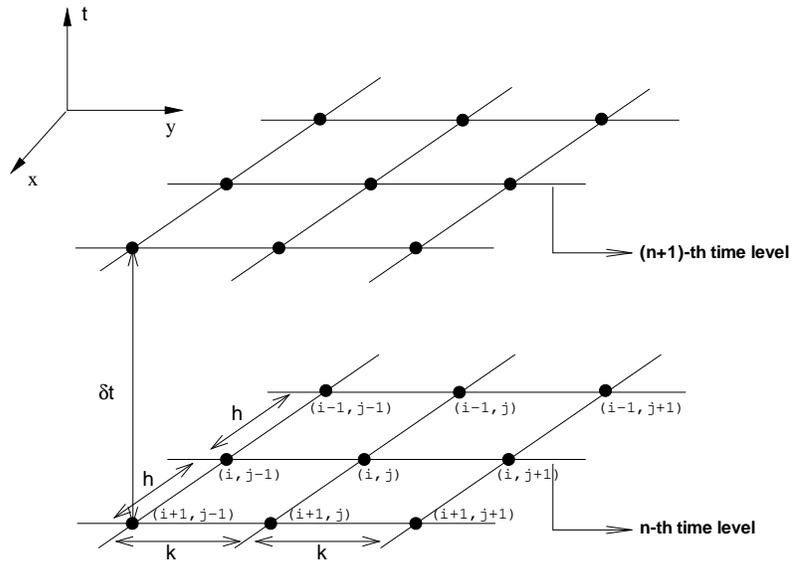}
\vspace{-1.5cm}
\caption{\sl Computational stencil: the used nodes are denoted by ``$\bullet$".}
    \label{fig:stencil}
\end{center}
\end{figure}

The HOC approach developed for the steady case can be extended directly to the unsteady case by simply replacing $s$ by $s-\partial_t \phi$ in Eq. (\ref{3}). At grid point $X_{i,j}$ and time $t$, the semi-discrete fourth order scheme for the parabolic equation with variable coefficients will be
\begin{eqnarray}\label{11}
\partial_t\phi_{i,j}(t)+A_{h,k}\phi_{i,j}(t)=s_{i,j}(t)
\end{eqnarray}
Clearly any time integrator can be used in Eq. (\ref{11}). Introducing weighted time average parameter $\iota$ such that $t_{\iota} = (1 - \iota)t^{(n)}_{\iota} + \iota t^{(n+1)}_{\iota}$ for $0 \leqslant \iota \leqslant 1$, where $(n)$ denote the $n$-th time level, we obtain a family of integrators; for example, forward Euler for $\iota=0$, backward Euler for $\iota=1$ and Crank-Nicholson for $\iota=0.5$. The resulting fully discrete difference scheme for grid point $(i, j)$ at time level $(n)$ then becomes
\begin{eqnarray}\label{12}
[1+\iota\delta tA_{h,k}]\phi^{(n+1)}_{i,j}=[1-(1 - \iota)\delta tA_{h,k}]\phi^{(n)}_{i,j}+\iota\delta t s^{(n+1)}_{i,j}+(1-\iota)\delta t s^{(n)}_{i,j}.
\end{eqnarray}
The accuracy of the scheme is $O(h^4, k^4, h^2k^2, \delta t^s)$, with $s\le 2$. The second order accuracy in time is obtained for $\iota=0.5$. On expansion the Eq. (\ref{12}) reduces to
\begin{eqnarray}\label{13}
&&[1-\iota\delta t(2\alpha^{(n+1)}_{1_{i,j}}\delta^2_x+2\alpha^{(n+1)}_{2_{i,j}}\delta^2_y-\beta^{(n+1)}_{_{i,j}}\delta_{x}\delta_{y}-d^{(n+1)}_{_{i,j}})]\phi^{(n+1)}_{_{i,j}}\nonumber\\
&&=[1+(1-\iota)\delta t(2\alpha^{(n)}_{1_{i,j}}\delta^2_x+2\alpha^{(n)}_{2_{i,j}}\delta^2_y-\beta^{(n)}_{_{i,j}}\delta_{x}\delta_{y}-d^{(n)}_{_{i,j}})]\phi^{(n)}_{_{i,j}}\nonumber\\
&&-(1-\iota)\delta t[(\alpha^{(n)}_{1_{i,j}}\delta_x-\beta^{(n)}_{_{i,j}}\delta_y+c^{(n)}_{1_{i,j}})\phi^{(n)}_{x_{i,j}}+(\alpha^{(n)}_{2_{i,j}}\delta_y-\beta^{(n)}_{_{i,j}}\delta_x+c^{(n)}_{2_{i,j}})\phi^{(n)}_{y_{i,j}}-s^{(n)}_{i,j}]\nonumber\\
&&-\iota\delta t[(\alpha^{(n+1)}_{1_{i,j}}\delta_x-\beta^{(n+1)}_{_{i,j}}\delta_y+c^{(n+1)}_{1_{i,j}})\phi^{(n+1)}_{x_{i,j}}\nonumber\\
&&+(\alpha^{(n+1)}_{2_{i,j}}\delta_y-\beta^{(n+1)}_{_{i,j}}\delta_x+c^{(n+1)}_{2_{i,j}})\phi^{(n+1)}_{y_{i,j}}-s^{(n+1)}_{i,j}].
\end{eqnarray}
It is clear that for all values of $\iota$, except $\iota=0$ and $1$, the difference stencil requires nine points in both $(n)$-th and $(n + 1)$-th time levels resulting in what may be called a $(9,9)$ scheme. The compact stencil emerging in this way has been illustrated in figure \ref{fig:stencil}. $(9,1)$ and $(1,9)$ schemes are obtained for $\iota=0$ and $\iota=1$, respectively. The algorithm given by Sen \cite{sen_13} can be used to solve the algebraic system associated with the finite difference approximation given in Eq. (\ref{13}).

\subsection{Stability analysis for constant coefficients}
We carry out the stability analysis for the scheme given by the equation Eq. (\ref{12}) by using von Neumann method. This is applicable in the special case when the coefficients $\alpha_1$, $\alpha_2$, $\beta$, $c_1$, $c_2$ and $d$ are constant. We consider a solution to the difference equation (\ref{13}) to be $$\phi^{(n)}_{i,j}=b^{(n)}e^{I\theta_{x}i}e^{I\theta_{y}j}$$ where $b^{(n)}$ is the amplitude at time level $n$, and $\theta_{x}=2\pi h/\Lambda_1$, $\theta_{y}=2\pi k/\Lambda_2$ are the phase angles with wavelengths $\Lambda_1$, $\Lambda_2$ respectively. Also periodic boundary conditions are assumed. The amplification factor for all $\theta_{x}, \theta_{y}\in[0,2\pi[$ is defined as $|G|:=|b^{(n+1)}/b^{(n)}|$.
\newtheorem{thm}{Theorem}
\begin{thm}[Stability]\label{theorem1}
The finite difference scheme (\ref{12}) with constant coefficients is stable, in the von Neumann sense, for $d\ge0$ if $\;0.5\leqslant\iota\leqslant 1$ and for $d<0$ the scheme is stable under sufficient condition $\displaystyle \delta t<-\frac{1}{\iota d}$.
\end{thm}
\newproof{pf}{Proof}
\begin{pf}
In order to prove the above theorem we will require the following two lemmas.
\newtheorem{lem}{Lemma}
\begin{lem}\label{lemma1}
$\alpha_1A^2+\alpha_2B^2+\beta ABr>0$ whenever $|r|\le1$.
\end{lem}
\emph{Proof of the Lemma 1:} This is clearly true if $\beta ABr>0$; otherwise since $\beta^2<4\alpha_1\alpha_2$ we get
\begin{eqnarray}
&&|\beta|<2\sqrt{\alpha_1\alpha_2}\nonumber\\
&\Rightarrow&|\beta ABr|<2\sqrt{\alpha_1\alpha_2}|AB|\nonumber\\
&\Rightarrow&-2\sqrt{\alpha_1\alpha_2}|AB|<\beta ABr\nonumber\\
&\Rightarrow&\alpha_1A^2+\alpha_2B^2-2\sqrt{\alpha_1\alpha_2}|AB|<\alpha_1A^2+\alpha_2B^2+\beta ABr\nonumber\\
&\Rightarrow&(\alpha_1|A|-\alpha_2|B|)^2<\alpha_1A^2+\alpha_2B^2+\beta ABr\nonumber\\
&\Rightarrow&0<\alpha_1A^2+\alpha_2B^2+\beta ABr.\nonumber
\end{eqnarray}
\begin{lem}\label{lemma2}
The function
\begin{eqnarray}
\mathfrak{R}(\theta_x,\theta_y)&=&\frac{2(8+\cos\theta_x+\cos\theta_y-\cos\theta_x\cos\theta_y)\cos(\theta_x/2)\cos(\theta_y/2)}
{\sqrt{(5+\cos\theta_x)(2+\cos\theta_x)(5+\cos\theta_y)(2+\cos\theta_y)}}\nonumber
\end{eqnarray}
is a periodic function of period $4\pi$ with respect to both the variables. This function defined over $(-\pi,3\pi)\times(-\pi,3\pi)$ attains (i) global maxima at $(0,0)$ and $(2\pi,2\pi)$ and its maximum value is 1, (ii) global minima at $(0,2\pi)$ and $(2\pi,0)$ and its minimum value is -1.
\end{lem}
We note that Lemma \ref{lemma2} can be arrived at, as a problem of maxima-minima.

Now we proceed to prove the stability theorem.

From equations (\ref{9}) and (\ref{10}) we see that
\begin{eqnarray}\label{14}
\phi_{{x}_{i,j}}^{(n)}=I\frac{3\sin\theta_{x}}{h(2+\cos\theta_{x})}\phi^{(n)}_{i,j},
\end{eqnarray}
\begin{eqnarray}\label{15}
\phi_{{y}_{i,j}}^{(n)}=I\frac{3\sin\theta_{y}}{k(2+\cos\theta_{y})}\phi^{(n)}_{i,j},
\end{eqnarray}
which upon substitution into the Eq. (\ref{8}) yields
\begin{eqnarray}\label{16}
A_{h,k}\phi^{(n)}_{i,j}=F(\alpha_1,\alpha_2,\beta,c_1,c_2,d,h,k,\theta_x,\theta_y)\phi^{(n)}_{i,j}
\end{eqnarray}
where
\begin{eqnarray}
&&F(\alpha_1,\alpha_2,\beta,c_1,c_2,d,h,k,\theta_x,\theta_y)\nonumber\\
\equiv && \bigg(2\alpha_1\frac{(2-2\cos\theta_x)}{h^2}+2\alpha_2\frac{(2-2\cos\theta_y)}{k^2}-\beta\frac{\sin\theta_x}{h}\frac{\sin\theta_y}{k}+d\bigg)\nonumber\\
&&+\bigg(-\alpha_1\frac{\sin\theta_x}{h}+\beta\frac{\sin\theta_y}{k}\bigg)\frac{3\sin\theta_x}{h(2+\cos\theta_x)}\nonumber\\
&&+\bigg(-\alpha_2\frac{\sin\theta_y}{k}+\beta\frac{\sin\theta_x}{h}\bigg)\frac{3\sin\theta_y}{k(2+\cos\theta_y)}\nonumber\\
&&+I\bigg(c_1\frac{3\sin\theta_x}{h(2+\cos\theta_x)}+c_2\frac{3\sin\theta_y}{k(2+\cos\theta_y)}\bigg)\nonumber\\
\equiv && F_{R}+d+IF_{I}.\nonumber
\end{eqnarray}
Thus the amplification factor of the scheme given in Eq. (\ref{12}) is then
\begin{eqnarray}\label{17}
&&|G|=\bigg|\frac{1-(1-\iota)\delta t F}{1+\iota\delta tF}\bigg|\nonumber\\
&\Rightarrow&|G|^2=\frac{(1-(1-\iota)\delta t (F_R+d))^2+(1-\iota)^2\delta t^2 F_I^2}{(1+\iota\delta t(F_R+d))^2+\iota^2\delta t^2F_I^2}.
\end{eqnarray}

Now for $d\ge0$ the scheme is unconditionally stable
\begin{eqnarray}
&&\rm{if}\;\;\;-2(F_R+d)+((F_R+d)^2+F_I^2)\delta t(1-2\iota)\le 0\nonumber\\
&i.e.&\rm{if}\;\;\;F_R\ge0\;\;\;\rm{for}\;\;\;(1-2\iota)\le0\nonumber\\
&i.e.&\rm{if}\;\;\;F_R\ge0\;\;\;\rm{for}\;\;\;0.5\le \iota\nonumber
\end{eqnarray}
Here,
\begin{eqnarray}
F_R&=&\frac{\alpha_1}{h^2}\bigg[4(1-\cos\theta_x)-\frac{3\sin^2\theta_x}{2+\cos\theta_x}\bigg]+\frac{\alpha_2}{k^2}\bigg[4(1-\cos\theta_y)-\frac{3\sin^2\theta_y}{2+\cos\theta_y}\bigg]\nonumber\\
&&+\frac{\beta}{hk}\sin\theta_x\sin\theta_y\bigg[\frac{3}{2+\cos\theta_x}+\frac{3}{2+\cos\theta_y}-1\bigg]\nonumber\\
&=&\frac{\alpha_1}{h^2}\bigg[\frac{\sin^2\theta_x-4\cos\theta_x+4}{2+\cos\theta_x}\bigg]+\frac{\alpha_2}{k^2}\bigg[\frac{\sin^2\theta_y-4\cos\theta_y+4}{2+\cos\theta_y}\bigg]\nonumber\\
&&+\frac{\beta}{hk}\sin\theta_x\sin\theta_y\bigg[\frac{8+\cos\theta_x+\cos\theta_y-\cos\theta_x\cos\theta_y}{(2+\cos\theta_x)(2+\cos\theta_y)}\bigg]\nonumber\\
&=&2\frac{\alpha_1}{h^2}\sin^2(\theta_x/2)\frac{5+\cos\theta_x}{2+\cos\theta_x}+2\frac{\alpha_2}{k^2}\sin^2(\theta_y/2)\frac{5+\cos\theta_y}{2+\cos\theta_y}\nonumber\\
&&+\frac{\beta}{hk}\sin\theta_x\sin\theta_y\bigg[\frac{8+\cos\theta_x+\cos\theta_y-\cos\theta_x\cos\theta_y}{(2+\cos\theta_x)(2+\cos\theta_y)}\bigg]\nonumber\\
&=&2\frac{\alpha_1}{h^2}\sin^2(\theta_x/2)\frac{5+\cos\theta_x}{2+\cos\theta_x}+2\frac{\alpha_2}{k^2}\sin^2(\theta_y/2)\frac{5+\cos\theta_y}{2+\cos\theta_y}\nonumber\\
&&+2\frac{\beta}{hk}\sin(\theta_x/2)\sin(\theta_y/2)\sqrt{\frac{5+\cos\theta_x}{2+\cos\theta_x}}\sqrt{\frac{5+\cos\theta_y}{2+\cos\theta_y}} \mathfrak{R}.\nonumber
\end{eqnarray}
By Lemma \ref{lemma2}, we see that $|\mathfrak{R}|\le1$. Hence by Lemma \ref{lemma1}, $F_R\ge 0$.

Next, $d<0$ correspond to a growth term in the governing equation (\ref{1}) and the stability condition requires $|G|<1+K\delta t$ with $K$ independent of $\delta t$ \cite{smith_book_86}. From (\ref{17}) we see that this is obviously true for $\displaystyle 1+\iota\delta t d>0\Rightarrow\delta t\le-\frac{1}{\iota d}$.

This completes the proof.
\end{pf}
\subsection{Extension to the 2D incompressible Navier-Stokes equations}
In this section we discusses applicability and effectiveness of the newly developed scheme in tackling coupled non-linear system of Navier-Stokes (N-S) equations in geometries beyond rectangular. It is worthwhile to point out that this fourth order compact scheme can be also used as a basis for discretizing unsteady incompressible N-S equations. One of the goals of this manuscript is to design schemes suited to simulate the time development of 2D viscous flows of varied physical complexities in different geometrical settings and accurately report fluid-embedded body interaction. Note that to resolve complex flow
phenomena accurately, we need to develop HOC schemes for body fitted coordinate system. It must be mentioned here that the idea of constructing a rectangular grid on irregular domains has its inherent limitations; furthermore accurate imposition of boundary condition becomes an issue. This is where the motivation to develop HOC schemes on nonuniform grids comes from. Although the primitive variable form of N-S equations can be discretized we choose to work with stream function-vorticity $(\psi-\omega)$ formulation as it decouples pressure field from velocity calculation. The system of equations presented below can be found elsewhere \cite{pan_kal_dal_07,sen_kal_gup_13} but we reproduce them for completeness.

In rectangular cartesian coordinates the $\psi-\omega$ formulation is given by
\begin{equation}\label{400}
-\psi_{xx}-\psi_{yy}=\omega,
\end{equation}
\begin{equation}\label{500}
\partial_t\omega-\frac{1}{Re}\big(\omega_{xx}+\omega_{yy}\big)+\big(u \omega_x + v \omega_y \big)=0.
\end{equation}
Here $\displaystyle (u,v)=\nabla^{\perp}\psi=\left (\psi_y, -\psi_x\right )$ is the velocity field and $\omega=-\nabla^2\psi$ is the vorticity. Also $Re=\frac{UL}{\nu}$ is the Reynolds number based on characteristic length $L$ and characteristic velocity $U$ of the flow. We consider non-singular mapping
\begin{equation}\label{600}
x=x(\xi,\eta),\;\; y=y(\xi,\eta).
\end{equation}
to transform the physical $(x, y)$ plane into a computational $(\xi, \eta)$ plane. Under this transformation equations (\ref{400}) and (\ref{500}) become
\begin{equation}\label{700}
-\widetilde{a}_{1}\psi_{\xi\xi}-\widetilde{e}_{1}\psi_{\xi\eta}
-\widetilde{b}_{1}\psi_{\eta\eta}+\widetilde{c}_{1}\psi_{\xi}
+\widetilde{d}_{1}\psi_{\eta}=\widetilde{f}_{1},
\end{equation}
\begin{equation}\label{800}
\partial_t\omega-\widetilde{a}_{2}\omega_{\xi\xi}-\widetilde{e}_{2}\omega_{\xi\eta}
-\widetilde{b}_{2}\omega_{\eta\eta}+\widetilde{c}_{2}\omega_{\xi}
+\widetilde{d}_{2}\omega_{\eta}=0,
\end{equation}
with
the velocity components expressed as
\begin{eqnarray}\label{900}
u=\frac{1}{J}\big(\psi_{\eta}x_{\xi}-\psi_{\xi}x_{\eta}\big),
\end{eqnarray}
\begin{eqnarray}\label{1000}
v=\frac{1}{J}\big(-\psi_{\xi}y_{\eta}+\psi_{\eta}y_{\xi}\big).
\end{eqnarray}
Details of the coefficients appearing in the above relations are available in \cite{pan_kal_dal_07,sen_kal_gup_13}. Although Pandit \emph{et al.} in their work \cite{pan_kal_dal_07} arrived at the above system of equations but ultimately restricted themselves to the development of HOC schemes for only conformably  transformed equations. Note that in case of conformal transformations coefficient of the mixed derivative term in the above equations will vanish. Here in this work we are not restricted by such choices of transformations and hence hold much wider appeal. To the best of our knowledge no fourth order compact formulation to tackle the above generalized N-S system is available in the literature. Both the equations (\ref{700}) and (\ref{800}) satisfy positive definiteness of associated diffusion matrix and are perfectly amenable to the discretization procedures detailed in this manuscript. The steady formulation given in Eq. (\ref{7}) is to be used for Eq. (\ref{700}) where as the implicit scheme presented in Eq. (\ref{12}) is to be used for Eq. (\ref{800}). In this manuscript we restrict ourselves to $\iota=0.5$ to have second order time accuracy. Since the gradients of the flow variables are already available at all grid points, necessity to compute $(u,v)$ separately does not arise. For the coupled nonlinear PDEs discussed above an outer-inner iteration procedure can be adopted. Outer one is the temporal iteration and in the inner iteration transient and steady equations are repeatedly solved till convergence is achieved.
\section{Numerical Examples}
\subsection{Problem 1: Analytical test case}
We consider the two-dimensional parabolic problem
\begin{eqnarray}\label{30}
\frac{\partial u}{\partial t}-\frac{\partial^2 u}{\partial x^2}+(1-x)(1-y)e^{x+y}\frac{\partial^2 u}{\partial x\partial y}-\frac{\partial^2 u}{\partial y^2}+10x(1-y)\frac{\partial u}{\partial x}-10y\frac{\partial u}{\partial y}=f(x,y,t)
\end{eqnarray}
defined over the domain $[0,1]^2\times(0,1)$. Following Karaa \cite{kar_07} we use this problem to test the proposed compact scheme and to confirm the theoretical fourth order accuracy in space and second order accuracy in time. In the Eq. (\ref{30}) the forcing term is selected to satisfy exact solution $u(x,y,t)=e^{-\pi t}(x^2-y^2)\cosh(x+y)$. The initial and Dirichlet boundary conditions are also taken from this exact solution. We consider three different mesh sizes $11\times11$, $21\times21$ and $41\times41$ and compute errors with respect to the exact solution using $L_1$, $L_2$ and $L_{\infty}$ norms. The computed results at time $t=0.25$ and $t=0.5$ have been presented in table \ref{table:P1_order_spatial}. A spatial order of accuracy close to four, for all cases, can be seen in this table. We estimate order of accuracy using the formula $\ln_2(Err1/Err2)$ where $Err1$ and $Err2$ are the errors with the grid sizes $h$ and $h/2$ respectively. Temporal second order accuracy of the scheme is verified in the table \ref{table:P1_order_temporal} where we compute using three different time steps $\delta t=0.01$, $\delta t=0.005$ and $\delta t=0.0025$ keeping $h=k=\tfrac{1}{30}$ unchanged.

\subsection{Problem 2: Boundary layer in non-rectangular domain}
This problem has been considered check the effectiveness of the spatial discretization in resolving boundary layers in non-rectangular domains. Here we tackle an elliptic problem by using the nine point scheme presented in equation (\ref{13}) with $\iota=0.5$. The domain of the problem is given by $$\Omega=\bigg\{(x,y):0\le x\le 1,0\le y\le \frac{1}{1-0.3\sin(6x)}\bigg\}.$$ The elliptic differential equation is
\begin{eqnarray}\label{40}
-\epsilon\bigg(\frac{\partial^2\psi}{\partial x^2}+\frac{\partial^2\psi}{\partial y^2}\bigg)-\frac{1.8\cos(6x)}{(1-0.3\sin(6x))(2-0.3\sin(6x))}\frac{\partial\psi}{\partial x}+\frac{1}{1+y}\frac{\partial\psi}{\partial y}=f,
\end{eqnarray}
which admits a boundary layer along the top curved surface of the domain $\Omega$. The analytical solution of the problem is $$\psi=e^{y-x}+(1+y)^{1+\tfrac{1}{\epsilon}}\bigg[\frac{1-0.3\sin(6x)}{2-0.3\sin(6x)}\bigg]^{\tfrac{1}{\epsilon}}.$$

\begin{figure}[!h]
\begin{center}
\includegraphics[width=2.5in]{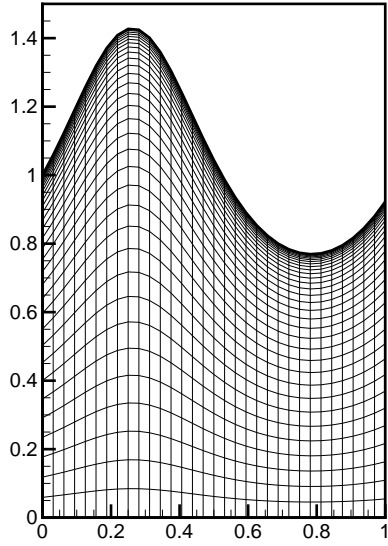}
\vspace{-1cm}
\caption{\sl Problem 2: $32\times32$ grid on the physical domain.}
    \label{fig:P2_grid}
\end{center}
\end{figure}

\begin{figure}[!h]
\begin{minipage}[b]{.6\linewidth}\hspace{-1cm}
\psfig{file=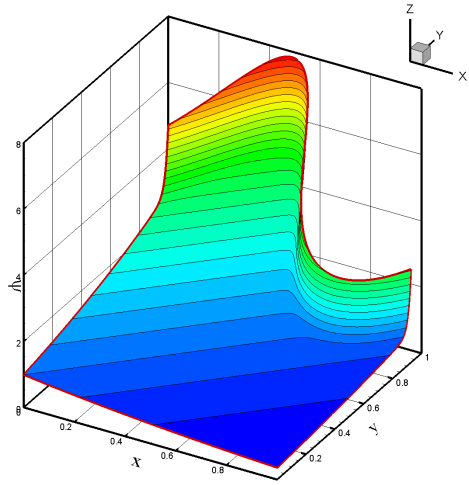,width=0.8\linewidth,angle=-90}
 \\(a)
\end{minipage}
\begin{minipage}[b]{.6\linewidth}\hspace{-1cm}
\psfig{file=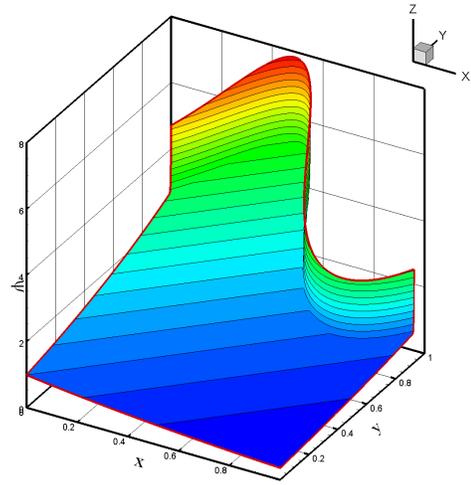,width=0.8\linewidth,angle=-90}
 \\(b)
\end{minipage}
\begin{center}
\caption{{\sl Problem 2: Numerical solution for (a) $\epsilon=0.01$, (b) $\epsilon=0.001$.} }
    \label{fig:P2_solution}
\end{center}
\end{figure}

A nonuniform grid along $y$-axis with clustering at the top and uniform grid along $x$-axis has been used following the transformation $$x=\xi,\;\;\;y=\frac{\eta+\tfrac{\lambda}{\pi}\sin(\pi\eta)}{1-0.3\sin(6\xi)}$$ from $(x,y)$ plane to $(\xi,\eta)$ plane. Here $\lambda$ is a parameter controlling the stretching along the top curved boundary and was considered to be 0.9 in our computations. With this transformation the equation (\ref{40}) takes the form presented in Eq. (\ref{13}). The physical domain $\Omega$ along with a $32\times32$ grids have been shown in figure \ref{fig:P2_grid}. We test our scheme using two small values of $\epsilon=0.01$ and $0.001$ using a fine grid of $257\times257$. The numerical results thus produced have been presented in figure \ref{fig:P2_solution}(a) and (b) respectively. The results clearly show that our computation has been able to capture the boundary layers in both the cases. For $\epsilon=0.01$ maximum error $\max||\textup{Num. Sol.}-\textup{Ana. Sol.}||_{\infty}\le10^{-8}$ where as for $\epsilon=0.001$ $\max||\textup{Num. Sol.}-\textup{Ana. Sol.}||_{\infty}\le10^{-7}$.
\subsection{Problem 3: Flow past circular cylinder}
\begin{figure}[!h]
\begin{minipage}[b]{.6\linewidth}\hspace{-1cm}
\psfig{file=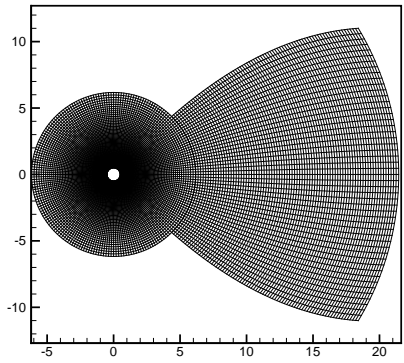,width=0.8\linewidth,angle=0}
 \\(a)
\end{minipage}
\begin{minipage}[b]{.6\linewidth}\hspace{-1cm}
\psfig{file=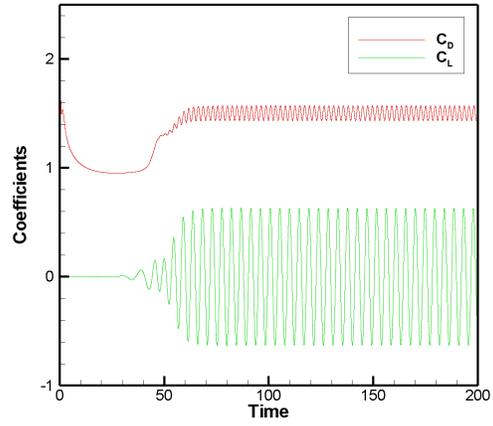,width=0.8\linewidth,angle=-90}
 \\(b)
\end{minipage}
\begin{center}
\caption{{\sl Problem 3: (a) Grid used for computing flow field. (b) Variation of drag and lift coefficients for $Re=300$.} }
    \label{fig:P3_grid_dl}
\end{center}
\end{figure}

\begin{figure}[!h]
\begin{center}
\includegraphics[width=4.0in,angle=-90]{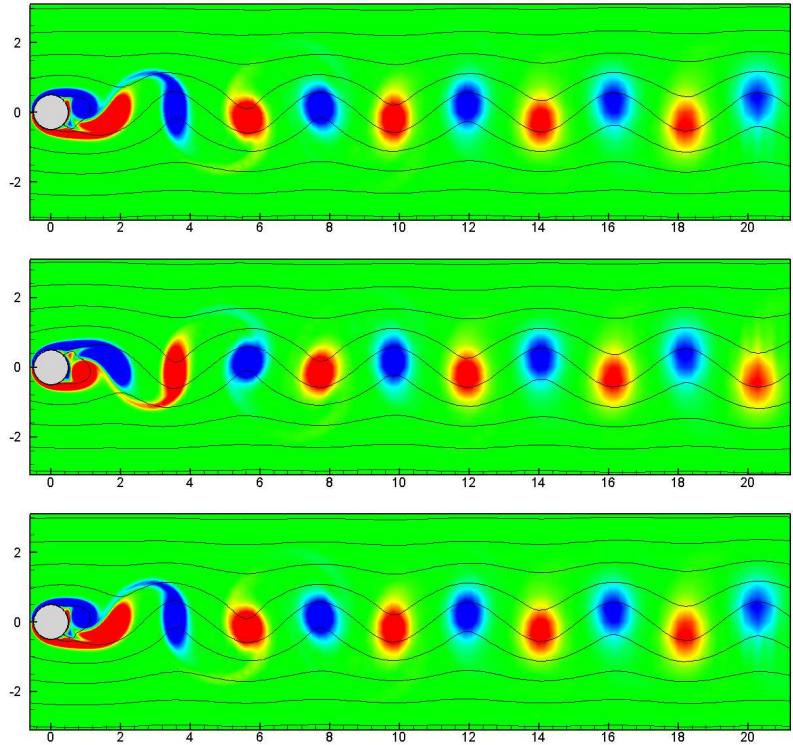}
\vspace{0cm}
\caption{\sl Problem 3: The streamlines and vorticity contours for three
successive instants of time (a) $t=T$, (b) $t=T+T_0/2$, (c) $t=T+T_0$.}
    \label{fig:P3_stream_vorti}
\end{center}
\end{figure}

Finally we intent to use the scheme developed here to solve Navier-Stokes equations and simulate flow around circular cylinder. This problem serves as an ideal test case for flow past bluff bodies on non-rectangular domains and of late, has quite often being used to validate schemes \cite{le_kho_per_06,ber_fal_08,wan_fan_cen_09,sen_sum_sin_10,pen_mit_sau_10} for their accuracy and effectiveness.  We assume the cylinder to be of unit diameter placed in an infinite domain and compute flow field for $Re=300$ by solving the coupled system of Eqs. (\ref{700}) and (\ref{800}). A multi-block CFD grid has been used for computing flow field. Precisely two blocks of grids have been used. Around the cylinder a circular grid of size $81\times201$ has been generated by using the transformation $x=\tfrac{1}{2}e^{(\pi\xi)}\cos(\pi\eta), y=\tfrac{1}{2}e^{(\pi\xi)}\sin(\pi\eta).$ Downstream of the flow, to capture the process of vortex shedding, another block of grids of size $81\times51$ have been numerically generated. The two sets of blocks have been point matched at the abutting boundary. The grid thus generated has been shown in the figure \ref{fig:P3_grid_dl}(a). The boundary conditions for the $\psi$ and its derivatives have been derived following the work of Kalita and Sen \cite{kal_sen_12-I}. Where as for $\omega$ we follow the boundary conditions given by Kalita and Ray \cite{kal_ray_09}. Tian \emph{et al.} \cite{tia_lia_yu_11} in their work have detailed various biased finite difference schemes for estimating first derivatives at the boundary. These approximations can be used here as well to approximate derivatives of vorticity field at the boundaries. We computing using $\delta t=0.005$. In figure \ref{fig:P3_grid_dl}(b), we show the time histories of the drag and lift coefficients which establishes the eventual periodic nature of the flow. In table \ref{table:P3_com} we compare Strouhal numbers, drag and lift coefficients for the flow with the values available in the literature. A good comparison can be seen. We present the temporal evolution of streamlines and vorticity contours over one complete vortex shedding cycle of duration $T=4.64$ in figure \ref{fig:P3_stream_vorti}. The evolution of von K$\acute{\rm a}$rm$\acute{\rm a}$n vortex street can be clearly seen in these figures. Two eddies are shed just behind the cylinder within each period. Figures \ref{fig:P3_stream_vorti}(a) and \ref{fig:P3_stream_vorti}(b) are half a vortex-shedding cycle apart and correspond to peak value of lift coefficient. Figure \ref{fig:P3_stream_vorti}(b) is a mirror image of figures \ref{fig:P3_stream_vorti}(a) and \ref{fig:P3_stream_vorti}(c). The crests and troughs of the sinuous waves in the streamlines reflect the alternatively positive and negative vorticities of the eddies.

\section{Conclusion}
In this paper, we present a higher order compact numerical method for 2D parabolic problems with mixed derivatives. Apart from solving convection-diffusion equation the formulation has been used to simulate unsteady two-dimensional Navier- Stokes equations. The robustness of the scheme is highlighted when it captures not only the boundary layer of a test case but also the process of vortex shedding for incompressible viscous flows on geometries beyond rectangular. Because of its high computational efficiency, our scheme has a good potential for efficient application. To the best of our knowledge compact schemes for parabolic 2D convection-diffusion equations is not available and thus this work is an important addition to the high accuracy solution procedures available in literature.

\bibliographystyle{plain}

\clearpage
\begin{table}[!h]
\begin{center}
\caption{\sl {Problem 1: $L_1$-, $L_2$-, $L_{\infty}$- norm error and spatial order of convergence with $\delta t=h^2=k^2$}}
{\begin{tabular}{ccccccc} \hline \hline
 Time  &      &$11\times11$ &Order &$21\times21$  &Order   &$41\times41$\\
\hline
$t=0.25$  &$L_1$  &2.338$\times10^{-6}$&3.88  &1.590$\times10^{-7}$&3.95  &1.034$\times10^{-8}$  \\
          &$L_2$  &2.972$\times10^{-6}$&3.91  &1.974$\times10^{-7}$&3.96  &1.269$\times10^{-8}$  \\
   &$L_{\infty}$  &5.823$\times10^{-6}$&3.94  &3.806$\times10^{-7}$&3.97  &2.428$\times10^{-8}$  \\
       &       &                    &      &                    &      &                    \\
$t=0.5$  &$L_1$  &1.067$\times10^{-6}$&3.88  &7.266$\times10^{-8}$&3.94  &4.724$\times10^{-9}$ \\
         &$L_2$  &1.357$\times10^{-6}$&3.91  &9.016$\times10^{-8}$&3.96  &5.798$\times10^{-9}$ \\
  &$L_{\infty}$  &2.657$\times10^{-6}$&3.93  &1.738$\times10^{-7}$&3.97  &1.109$\times10^{-8}$ \\
\hline \hline
\end{tabular}}
\label{table:P1_order_spatial}
\end{center}
\end{table}
\clearpage
\begin{table}[!h]
\begin{center}
\caption{\sl {Problem 1: $L_1$-, $L_2$-, $L_{\infty}$- norm error and temporal order of convergence with $h=k=\tfrac{1}{30}$}}
{\begin{tabular}{ccccccc} \hline \hline
 Time  &      &$\delta t=0.01$ &Order &$\delta t=0.005$  &Order   &$\delta t=0.0025$\\
\hline
$t=0.25$  &$L_1$  &8.050$\times10^{-7}$&1.95  &2.088$\times10^{-7}$&1.73  &6.287$\times10^{-8}$  \\
          &$L_2$  &1.054$\times10^{-6}$&1.94  &2.746$\times10^{-7}$&1.72  &8.356$\times10^{-8}$  \\
   &$L_{\infty}$  &2.406$\times10^{-6}$&1.93  &6.294$\times10^{-7}$&1.75  &1.872$\times10^{-7}$  \\
       &       &                    &      &                    &      &                    \\
$t=0.5$  &$L_1$  &3.671$\times10^{-7}$&1.95  &9.524$\times10^{-8}$&1.73  &2.869$\times10^{-8}$ \\
         &$L_2$  &4.806$\times10^{-7}$&1.94  &1.253$\times10^{-7}$&1.72  &3.814$\times10^{-8}$ \\
  &$L_{\infty}$  &1.097$\times10^{-6}$&1.93  &2.871$\times10^{-7}$&1.75  &8.542$\times10^{-8}$ \\
\hline \hline
\end{tabular}}
\label{table:P1_order_temporal}
\end{center}
\end{table}
\clearpage
\begin{table}[!h]
\begin{center}
\caption{\sl {Comparison of Strouhal numbers, drag and lift
coefficients of the periodic flow for  $300$.}}
\begin{tabular}{cccc}
\hline
Reference  &$St$&$C_D$ &$C_L$
\\\hline
Frank {\it et al.} \cite{fra_rod_sch_90}&0.205&1.32&$\pm$ 0.84 \\
Wang \emph{et al.} \cite{wan_fan_cen_09} &0.206&1.174$\pm$ 0.080&$\pm$ 0.90 \\
Sen \emph{et al.} \cite{sen_kal_gup_13} &0.209&1.401$\pm$0.068&$\pm$0.607\\
Present Study&0.215&1.503 $\pm$ 0.069 &$\pm$0.630 \\
\hline
\end{tabular}
\label{table:P3_com}
\end{center}
\end{table}
\clearpage

\clearpage

\clearpage
\clearpage

\clearpage

\clearpage
\end{document}